**Werner DePauli-Schimanovich**[1]

**Higher Cardinals are only a Convention, Ch64-SE11**

Abstract:
Zermelo's Axiom of Separation is:
Exist x: Forall y: (y in x <==> y in a & E(y)) with definite(E) for a parameter a.
Thoralf Skolem suggested to characterize the terminus "definite" by "the property E should be representable by a FOL formula". But that is trivial. "definite" must mean more. The author claims that "definite" means "in accordance with the theory of definitions of logic". In this case the theorem of Cantor is no longer a theorem, but a undecidable sentence, and has to be established explicitly as axiom. This is not done by the community, but it is made a silent assumption that we can drop the appendix "definite(E)" from the axiom of separation at all. But this is a convention (even when it is silent) and it is nothing else than an axiom.

Paper:
At first, Georg Cantor, the inventor of set theory, still believed that "in the night of infinity all cats are grey!" (i.e."at night, everything is the same color").[2] He knew that rational numbers are countable. But he could not demonstrate the cardinality of real numbers, and did not consider this question particularly important. Still, in his letter to Richard Dedekind on December 7$^{th}$ in 1873[3], Cantor provided proof for the uncountability of real numbers.

The whole real numbers R are equal in cardinality to the interval [0,1] and can be seen as infinite binary fractions between 0 and 1 instead of as infinite decimal fractions. These binary fractions in turn constitute the number of functions of the natural numbers ℕ in {0,1}, which is nothing other than the power set ℙ(.) of ℕ. With this, the basis for generalization was already set, i.e., for Cantor's theorem[4]: [For all a:] $|a| < |\mathbb{P}(a)|$.
The cardinality of a is smaller than the cardinality of its power set.

For proving the uncountability of real numbers between 0 and 1, Cantor constructs a "diagonal element" from the list of all real numbers between 0 and 1, which is presupposed to exist, and from that he builds an "anti-diagonal element" ADE that also would have to be a real number, but nevertheless cannot be included in the list because of its construction. Here Cantor's opponents brought a well-known counter-argument – that it is impossible to suddenly construct a new real number out of a closed list of all real numbers, but the new number is not a member of this closed list, in spite of the fact that it would have to be, since this is an implicit presupposition!

Naturally, this objection comes from the constructivists. Platonists don't give it a thought, since they see the list as an actual existing infinite set. This led to the wrong conclusion that the power set axiom should be rejected.[5] But the problem is the axiom of separation (i.e., the partial-set axiom) which is formulated by Ernst Zermelo only for "definite" properties. But what should that mean: a definite property?[6] The author is convinced that it can have only one meaning: "in accordance with the rules of the theory of definitions". And this means mainly that the definiendum (= result of the definition) is not allowed to appear in the definiens except in recursive definitions.

Thoralf Skolem gave another formulation to make the terminus "definite" precise: Definite properties are exactly those unary predicates which can be defined with parameters. Or: Definite(E) := ∃ sets $x_1, ...x_n$ and a wff $\varphi(z, x_1, ... x_n)$ such that "t(z) is valid ⇔ $\varphi(z, x_1, ... x_n)$". See Ebbinghaus [1994], page 30, or Skolem [1922]. But I would object this view. At least we have to add "… concerning the rules of definition", because Skolem's formulation alone constitutes only that the definite property E can be expressed as a predicate logical formula (of 1$^{st}$ order). But since in logic and in set theory ZF (codified in FOL) all properties can be expressed as predicate logical formulas, Skolem's formulation of "definite" is a tautology. Therefore it is dropped today in the axiom of separation and also replacement (i.e. image-set axiom). But if we drop it we say implicitly that "definite" has no meaning and that was definitely not Zermelo's intention of the axiom of separation.

---

[1] Dept. of DB&AI, Technical University Vienna, Favoritenstraße 9, A-1040 Wien, Austria/Europe. Werner.DePauli@gmail.com

[2] See Meschkowski [1973] page 16.
[3] Meschkowski [1973], page 17, considers this date to be the birthday of set theory.
[4] |a| is the cardinality (or power) of a, "<" means "…is of lower power than…", and ℙ is the power set (i.e. set of all subsets) of a given basic set. Sometimes when we want to say it explicitly that "<" means really smaller, we write also "≨".
[5] For an introduction to axiomatic set theory, see Suppes [1960] or Ebbinghaus [1994].
[6] Juliet Floyd and Akihiro Kanamori write in [2006], page 42 "But what exactly is a "definite" property? This was a central vagary … of the Zermelo Set Theory."

Either the reader agrees with me that "definite" should mean "according to the rules of definition" (what would be the correct way), or he/she must allow also properties E which are not formed according to the rules of definition. Of course we can do so. But this is something like a new axiom: A complete new situation that we need not to obey the rules of definition. In any case it is an agreement of the community of mathematicians that in this case of the axiom of separation (and also of replacement) we need not to form the formula E according the theory of definitions. This is no selfevident triviality, but a convention! Because the theory of definitions is part of logic and more evident than any mathematical formula or axiom.

Sometimes parameters which are not well defined concerning to the rules of the theory of definitions are also called "impredicative"[7] parameters.

Nevertheless, it would be correct to forbid impredicative parameters in the axiom of separation.

The above-mentioned generalization in set theory corresponds to Cantor's theorem concerning real numbers. Here, too, proof is derived indirectly: the hypothesis that there exists a function from a on $\mathbb{P}(a)$ leads to a contradiction. Here, too, an anti-diagonal element is defined.

$ADE(f) = \{x \in a \,|: x \notin f(x)\}$.

$ADE(f)$ is, according to the axiom of separation (as usually interpreted), a subset of a, and therefore an element of the power set $\mathbb{P}(a)$. It is immediately obvious, however, that for no $x_0 \in a$ can the function value $f(x_0) = \{x \in a \,|: x \notin f(x)\}$ exist, for then

$x_0 \in f(x_0) \iff x_0 \notin f(x_0)$ would be true: a contradiction[8]; therefore, there exists no such function. This is how the standard argument proceeds.

However, the fact that a second hypothesis was made for this proof was assiduously overlooked - that is, the axiom of separation was interpreted such that $ADE(f)$ should be a set in all cases. It bases on this contradiction, we can see that either no function exists, or $ADE(f)$ is not a set because $x \notin f(x)$ is not definite! It could be a properclass such as the Russell-class, which would actually be understandable, due to the use of impredicative parameters. The construction of f depends upon $ADE(f)$, and $ADE(f)$ presupposes a completed f. This produces a circle in the sequence of elements and is against the rules of definitions.

We have to be cautious: Of course we have always the relation

$$f : a \to b \equiv \{\langle x, f(x) \rangle \,|: x \in a \,\&\, f(x) \in b\} = f =$$
$$= \{\langle x, y \rangle \,|: \exists x \in a, \exists y \in b \,\&\, f(x) = y\} = f.$$

But defining f, b must be independent from f, because it's a free parameter of "f: a $\to$ b". This is not the case here because $ADE \in b$ is defined by $ADE := \{x \in a \,|: x \notin f(x)\}$. This is against the rules of definition because f is defined as reflexive relation and not recursively. Such a relation can be an axiom but is not allowed to form a definition. The definiendum f cannot appear in the definiens (of which ADE is a part) except in a recursive definition. Therefore ADE cannot even be defined correctly concerning to the rules of definition.

Furthermore ADE has nothing to do with the higher cardinality of the power set. Even the function
$f: \{\emptyset\} \to \{\{x \in \{\emptyset\}|: x \notin f(x)\}\}$
with one-elemented domain and one-elemented range does not exist. But such a function has to exist in any case: Either it is the identity 1 restricted to $\{\emptyset\} \times \{\emptyset\} := \{\{\emptyset\}\}$, or the function $\{\langle \emptyset, \{\emptyset\} \rangle\}$. If we do not want to consider this as an inconsistency within ZFC, we should interpret the terminus definite in my way. But if we interpret "definite" as "following the theory of definitions", we cannot prove any longer the "Theorem of Cantor". Its now an undecidable sentence like the one in Gödel's proof.[9]

What is the conclusion of that? We have to agree that the "Theorem of Cantor" is no theorem, but a convention between the mathematicians, which has to be added as axiom to ZFC. The best way to regulate that would be to assume the Generalized Continuums Hypotheses (GCH). Despite of the fact that some very counter-intuitive facts[10] follow from (GCH) and that outstanding logicians[11] had been therefore against (GCH), it makes the cardinal arithmetic much simpler. And what other reason should we have to reject it if even the higher cardinality of the powerset is a convention?

---

[7] The terminus "impredicative" has the first time be used by Henry Poincare but had never been precised.
[8] This contradiction has great similarities to Russell's Paradox. In fact it was Cantor's proof that delivered Bertrand Russell the idea for his antinomy: Set a := V and f := id and ADE becomes the Russell class $\mathbb{R}u$.
[9] See my books about Gödel in the literature.
10 From (GCH) follows e.g. the so-called Banach-Tarski paradoxon: you can decomposite the volume of a ball into pieces that form 2 congruent balls (i.e. dense multitudes or sets of points) where each of it have the same volume and diameter as the original one. See e.g. http://de.wikipedia.org/wiki/Banach-Tarski-Paradoxon.

The axiom of separation (as well as the axiom of replacement) ought in any case to be limited, because it inconveniently hinders the formation of sets due to its inherent "limitation of size" ideology[12], and allows for no universal sets.[13] For example, the basis set a should be limited to well-founded sets. This and two other small changes make it possible to add the existence of complements to such a restricted ZFC.[14] In addition, we could consider requiring the stratification[15] of separating property B. Thus:
(For all stratified wff B-Strat $\exists$ a B-Strat): $\exists$ x $\forall$ y : y $\in$ x $\Leftrightarrow$ y $\in$ a & B-Strat(y).

The proof of Cantor's theorem described above can only work in this way because the set-generating comprehension-operator {x: A(x)} is able to contain the well-formed formula A with an impredicative (free) parameter. If we required the stratification of partial formula B, that would not be possible. As we know, the use of parameters in mathematics is a perfectly typical way of doing business to which nobody objects. Nevertheless, it is in no way trivial that impredicative parameters may be used without limit also in the wff A in the Church Comprehension-Scheme: y $\in$ {x: A(x)} $\Leftrightarrow$ A(y), or in the scheme of separations.

In Zermelo-Fraenkel's set theory ZFC with an axiom of choice (AC), there is a famous model: the cumulative hierarchy. It is constructed bottom up starting with the empty set, from which the power set is formed again and again. All iterations are unified, from which once again the power set is formed, etc. As a whole, this results in the universal class V.
    In constructing this cumulative hierarchy, no parameters at all must be employed - here, this is unnecessary. And the problems of the anti-diagonal element do not appear, either. The ADE(f) depends on parameter f. If we then form $f(x_0)$ = ADE(f), however, we land in a *circulus vitiousus*: we first need f in order to form ADE. In set theory, functions are sets of ordered pairs with the "functional property". This results in:
$<x_0, ADE(f)> \in$ f, i.e. f is not well-founded.
    If we assume also the axiom of foundation and allow impredicative properties, this fact would already be sufficient for the higher power of the power set. We need not to show in addition that there is no function between the set and its power set.

We can accept the "Theorem of Cantor" as an axiom, even such a self-reference as in ADE by cancelling the "definiteness"-condition in Zermelo's axiom of separation (or Fraenkel's axiom of replacement) with confidence, since we know from experience that it does not lead to any contradiction, but only to higher cardinalities. However, it must be clear to us that this is purely convention and does not constitute any general truth. That republic of experts known as mathematicians has quietly decided to always allow the use of parameters in set operators - even when they are impredicative and sometimes also indefinite and (as in the example of our function f) therefore create a self-reference.[16] Only through this convention exist differing higher cardinalities in the transfinite.

My intention here is not to cast the existence of higher cardinalities into doubt. However, it does not constitute a Platonic truth, but rather a convention. Since over 90% of mathematicians are Platonists (for whom the actual-infinite has a real existence), this issue never is questioned by them. But it is a mistake to act as if higher cardinalities were a God-given triviality! There is no question that we recognize this convention in mathematics. Nevertheless, there is no contradiction should we also recognize another convention. Physicists, cosmologists, quantum theoreticians, and other scientists (as well as economists) all have no need of higher cardinalities and therefore do need some option other than ZFC with impredicative parameters in terms of sets.[17]

---

[11] Also Gödel was against (GCH). See his paper "What is Cantor's Continuum Problem" in the Collected Works.
[12] See Hallet [1984].
[13] See Forster [1995]
[14] See DePauli-Schimanovich [2008a].
[15] A wff A is called stratified, if the variables can be indexed with natural numbers such that for every pair of variables left and right of an element-sign "$\in$" the index on the right is 1 higher than the one on the left side. E.g. "x7 $\in$ y8 &…& z6 $\in$ x7" is stratified, but "x $\in$ y &…& y $\in$ x" not.
[16] With self-references, we must be very careful, since they can easily give rise to paradoxes (e.g., Russell's Paradox).
[17] ZFC shall be here the extended system with Definition by Abstraction and the use of set-operator. See Suppes, p. 33.

For scientists, and in particular for physicists, a set theory with two transfinite cardinalities is sufficient - countability and continuum. Randall Holmes[18], for example, created a Pocket Set Theory that has only these two higher cardinalities.

If in the set operators in the set-constituting wffs $A_i$ of the axiom of replacement and separation are not allowed any impredicative parameters (but only, for example, stratified formulas), Cantor's theorem cannot be proven in the usual straightforward way.[19] Then we have a choice: either we use the Generalized Continuum Hypothesis GCH, thus fixing the higher cardinalities[20], or we choose a simple CH with countable and continuous sets (as in Pocket Set Theory), or we can choose an anti-hypothesis AH, which stipulates that no transfinite cardinalities may exist other than countable sets. We could even define the set $\omega_0$ of the natural numbers as a properclass and thus allow only finite sets, as the quantum mechanists would have it.

Since the concept "impredicative" is difficult to formalize, it is simpler to prohibit all parameters other than the basis set a, or to limit them suitably. It is sufficient to prohibit or restrict the parameters in the axiom of replacement (i.e., image-set axiom), since the axiom of separation can be logically derived from it. In this connection, I would like to also mention my maximal possible set theory NACT*[21], whose class-terms allow no parameters whatsoever, and whose sets are formed in a manner similar to cumulative hierarchy. The invalidity of Cantor's theorem is in this case not a *malheur*, but rather a *bonheur*, for now we are presented with the choice of how many transfinite cardinalities we wish to have.

---

[18] See Holmes [2005].

[19] In NF (= New Foundation) of Quine or NFUM (= New Foundation with Urelements and Measurable Ordinal-class) of Holmes [1998], Cantor's theorem can be proved with the help of some technical tricks. But this proof has the same problems we discuss here in this paper.

[20] This is possible due to the relative consistency of (GCH) with ZFC, which has been proven by Kurt Gödel. Concerning Gödel, see Casti/DePauli [2000], Weibel/Schimanovich [1986] and Buldt et al [2002].

[21] See DePauli-Schimanovich [2008c], [2002b], [2006b].